# An Efficient Finite Difference-based Implicit Solver for Phase-Field Equations with Spatially and Temporally Varying Parameters

Zirui Mao[†a], G. R. Liu[b], Michael J. Demkowicz[a]

[a]Department of Materials Science and Engineering,

Texas A&M University, College Station, Texas 77843, USA

[b]Department of Aerospace Engineering and Engineering Mechanics,

University of Cincinnati, Cincinnati, Ohio 45221, USA

**Abstract**: The phase field method is an effective tool for modeling microstructure evolution in materials. Many efficient implicit numerical solvers have been proposed for phase field simulations under uniform and time-invariant model parameters. We use Eyre's theorem to develop an unconditionally stable implicit solver for spatially non-uniform and time-varying model parameters. The accuracy, unconditional stability, and efficiency of the solver is validated against benchmarking examples. In its current form, the solver requires a uniform mesh and may only be applied to problems with periodic, Neumann, or mixed periodic/Neumann boundary conditions.

## 1. Introduction

The phase field method (PFM) [1], [2] is a mathematical model that can describe and track the evolution of microstructure in multi-phase materials by representing distributions of different phases as field variables. It has been successfully used in a wide range of microstructure evolution problems, e.g., solidification [3], grain growth [4], [5], phase separation [6]–[8], crack propagation [9], [10], and crystal growth [11][12]. It is readily applied to a variety of solids, including geological materials [13], [14], engineering composites [15], [16], shape memory alloys [17], [18], and Li-ion battery components [19]–[21].

Generally, phase field (PF) equations are a group of second-order or fourth-order partial differential equations (PDEs). The finite difference (FD) method [22] is one common way of solving PF boundary value problems. To enhance computational efficiency, stable implicit solvers allowing larger time steps have been developed [23]–[32].

[†] Correspondence to: Zirui Mao, Department of Materials Science and Engineering, Texas A&M University, College Station, Texas 77843, USA. maozu@mail.uc.edu

These methods have been developed to handle PF problems with fixed model parameters. Physically, this constraint corresponds to a time-invariant external state, e.g., a stable environment (temperature, pressure, chemical potential) or fixed processing conditions (mechanical load, electromagnetic fields). However, microstructure evolution often takes place under changing external conditions. Varying parameters cause a more challenging derivation of stability condition of numerical solvers because model parameters in different range lead to different stability criteria. Accordingly, there is a need for an efficient, stable solver for PF equations with varying parameters. We propose an unconditionally stable implicit FD solver based on the von Neumann numerical stability criterion [33] and Eyre's gradient stability theorem [34]. The proposed implicit scheme may be solved directly, yielding greater efficiency than iterative methods.

The remainder of this paper is organized as follows. Section 2 defines two model problems to solve in this study. Section 3 presents the proposed implicit FD solver and its direct solution. Section 4 investigates the accuracy, stability, and scalability of the proposed implicit FD solver when solving the two model problems defined in section 2. We conclude with a discussion of the merits and limitations of the proposed implicit FD solver in section 5. A theoretical deduction of the proposed implicit scheme is provided in the appendix.

## 2. Phase field equations

This paper will concentrate on two classical PF equations with scalar fields as model problems: the second order Allen-Cahn (A-C) equation and the fourth order Cahn-Hilliard (C-H) equation. The A-C equation [35], [36] has the form

$$\frac{\partial \phi}{\partial t} = -M\mu, \ \mu = \frac{\partial f(\phi)}{\partial \phi} - \gamma \nabla^2 \phi \tag{1}$$

where $\phi$ is a field variable called the 'order parameter' (a non-conserved quantity), $M$ is mobility, $f$ is a free energy function that depends on $\phi$, and $\gamma$ is the gradient energy coefficient. The Cahn-Hilliard (C-H) equation [37] has the form

$$\frac{\partial \phi}{\partial t} = \nabla \cdot (M \nabla \mu), \ \mu = \frac{\partial f(\phi)}{\partial \phi} - \gamma \nabla^2 \phi \tag{2}$$

where the order parameter $\phi$ is now a conserved quantity and the remaining terms have the same meaning as in the A-C equation.

The A-C and C-H equations are derived from the definition of gradient flow [38], where the evolution of $\phi$ follows the path of steepest descent of the total Ginzburg-Landau free energy, $F$, in domain $\Omega$ [37]:

$$F(\phi) \coloneqq \int_\Omega \left( f(\phi) + \frac{\gamma}{2} |\nabla \phi|^2 \right) dx \tag{3}$$

Hence, the proposed FD scheme must satisfy not only numerical stability, but must also guarantee the monotonous decay of *F* in the evolution process, which, in physics, is called 'gradient stability.'

In both A-C and C-H equations, the second term in the definition of the chemical potential, $\mu$, reflects the gradient energy stored in the interface between adjacent phases, where the order parameter varies rapidly. The gradient energy coefficient $\gamma$ may be viewed as a constitutive parameter (material property). The first term in $\mu$ reflects the uniform (gradient-independent) part of the free energy density, *f*. This term depends directly on $\phi$ as well as any number of parameters describing the external environment, such as temperature or various fields.

In this study, we set

$$f = \phi^4 + T(\pmb{x},t)\phi^2 + h(\pmb{x},t)\phi \tag{4}$$

where *T* and *h* are external parameters, each of which may vary with location as well as time. As illustrated in Figure 1.(a), *T* has a similar effect to temperature, because large *T* favors a free energy with a single minimum, corresponding to just one equilibrium value of $\phi$. By contrast, low *T* favors a free energy with two minima, i.e., two equilibrium $\phi$ values, leading to phase separation. Parameter *h* has a similar effect as a symmetry-breaking field, as non-zero *h* gives an asymmetric *f*, shown in Figure 1.(b), favoring the phase with lower *f*.

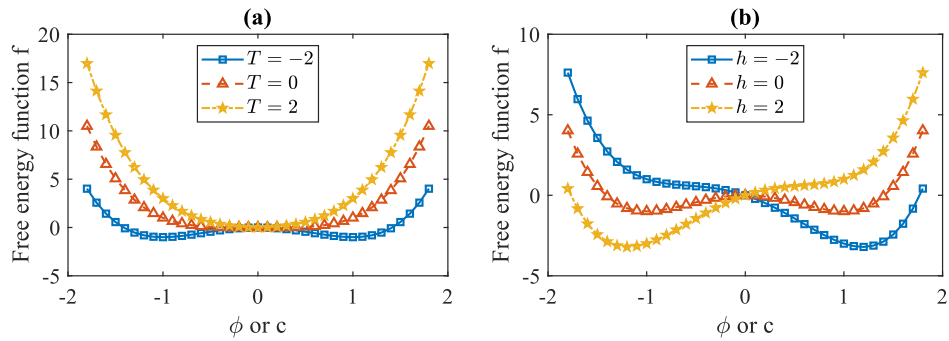

Figure 1: The free energy density in Eqn. 4 for selected values of *T* and *h*. In (a), *h* = 0. In (b), *T* = -2.

While we illustrate our solver on the free energy density function in Eqn. (4), it is readily applied to other forms of *f*.

## 3. Finite difference solvers

In this section, we first present the explicit FD solver for solving A-C and C-H equations with varying parameters. On this basis, we develop our unconditionally stable implicit FD solver as well as the associated direct solution to the proposed implicit FD schemes. Our presentation uses 2-D domains, but the method is easily extended to 3-D.

## 3.1 Explicit FD solver

The PF PDEs in equations (1)-(2) can be solved by approximating the Laplacian terms with the 2$^{nd}$ order accurate central finite-difference scheme:

$$\nabla^2 u_{i,j} = \frac{u_{i,j+1} - 2u_{i,j} + u_{i,j-1}}{\Delta x^2} + \frac{u_{i+1,j} - 2u_{i,j} + u_{i-1,j}}{\Delta y^2} + O(\Delta x^2 + \Delta y^2) \tag{5}$$

If assuming equal spacing *h* in both *x* and *y* directions, *i.e.*, $\Delta x = \Delta y = \Delta r$, the above equation becomes:

$$\nabla^2 u_{i,j} = \frac{u_{i+1,j} + u_{i-1,j} + u_{i,j+1} + u_{i,j-1} - 4u_{i,j}}{\Delta r^2} + O(\Delta r^2) \tag{6}$$

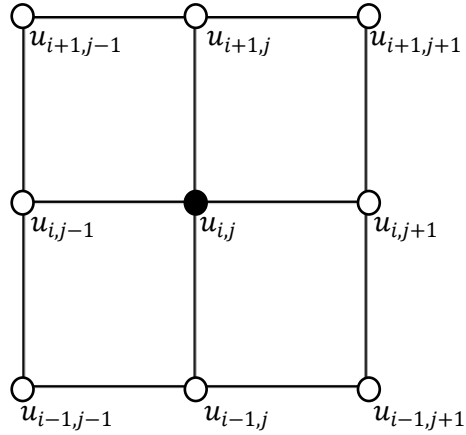

Figure 2 A schematic node ordering representation for 2-D finite difference scheme.

For the equally spaced nine grid points shown in Figure 2, and without considering boundary conditions, Eq. (6) can be rewritten in a matrix form as:

$$\nabla^2 \mathbf{u} = \frac{1}{\Delta r^2} \left( \begin{bmatrix} -2 & 1 & & & & 0 \\ 1 & -2 & 1 & & & \\ & 1 & \ddots & \ddots & & \\ & & \ddots & \ddots & 1 & \\ & & & 1 & -2 & 1 \\ 0 & & & & 1 & -2 \end{bmatrix} \mathbf{u} + \mathbf{u} \begin{bmatrix} -2 & 1 & & & & 0 \\ 1 & -2 & 1 & & & \\ & 1 & \ddots & \ddots & & \\ & & \ddots & \ddots & 1 & \\ & & & 1 & -2 & 1 \\ 0 & & & & 1 & -2 \end{bmatrix} \right) = \mathbf{A}\mathbf{u} + \mathbf{u}\mathbf{A} \tag{7}$$

where the coefficient matrix **A** is symmetric, tri-diagonal, and invertible. Therefore, the solution exists and is unique. Convergence is ensured by the Taylor expansion theory, because the grid used here is regular.

Consequently, the A-C equation in eq. (1) can be solved as the following time-evolution equation

$$\boldsymbol{\phi}^{t+1} = \boldsymbol{\phi}^t + \Delta t \left[ -M(f'_{\phi^t} - \gamma(\mathbf{A}\boldsymbol{\phi}^t + \boldsymbol{\phi}^t \mathbf{A})) \right]$$

and the C-H equation in eq. (2) can also be solved by

$$\boldsymbol{\phi}^{t+1} = \boldsymbol{\phi}^t + M\Delta t \, (\mathbf{A}\boldsymbol{\mu}^t + \boldsymbol{\mu}^t \mathbf{A}), \text{ with } \boldsymbol{\mu}^t = f'_{\phi^t} - \gamma(\mathbf{A}\boldsymbol{\phi}^t + \boldsymbol{\phi}^t \mathbf{A})$$

Boundary conditions (BCs) may be incorporated into the coefficient matrix **A**. For periodic, Neumann with zero-gradients, and symmetric BCs, **A** has the forms of $A_p$, $A_N$, and $A_S$, respectively:

$$\mathbf{A}_p = \frac{1}{\Delta r^2}\begin{bmatrix} -2 & 1 & 0 & \cdots & & 0 & 1 \\ 1 & -2 & 1 & 0 & & \cdots & 0 \\ 0 & 1 & \ddots & \ddots & & \ddots & 0 \\ \vdots & \ddots & & \ddots & \ddots & 1 & 0 \\ 0 & \cdots & & 0 & 1 & -2 & 1 \\ 1 & 0 & & \cdots & 0 & 1 & -2 \end{bmatrix} \qquad (8)$$

$$\mathbf{A}_N = \frac{1}{\Delta r^2}\begin{bmatrix} -1 & 1 & & & & & \\ 1 & -2 & 1 & & & 0 & \\ & 1 & \ddots & \ddots & & & \\ & & \ddots & \ddots & 1 & & \\ & & & 1 & -2 & 1 & \\ & 0 & & & & 1 & -1 \end{bmatrix} \qquad (9)$$

$$\mathbf{A}_s = \frac{1}{\Delta r^2}\begin{bmatrix} -2 & 2 & & & & & \\ 1 & -2 & 1 & & & 0 & \\ & 1 & \ddots & \ddots & & & \\ & & \ddots & \ddots & 1 & & \\ & & & 1 & -2 & 1 & \\ & 0 & & & & 2 & -2 \end{bmatrix} \qquad (10)$$

It is straightforward to apply different BCs in different directions. For example, $\nabla^2 \boldsymbol{u} = 1/\Delta r^2 (A_p \boldsymbol{u} + \boldsymbol{u} A_S)$ represents periodic BCs in the vertical direction and a symmetric BC in the horizontal direction.

The explicit FD solver is simple to implement, but the time step $\Delta t$ is subject to the von-Neumann stability condition. Consequently, the explicit solver is relatively inefficient, especially for high-resolution problems with a very small $\Delta r$. In 2-D, the stability criterion of the explicit FD solver for A-C and C-H equations has the following forms [39], respectively:

$$\Delta t_{A-C} \leq \frac{\Delta r^2}{4\gamma} \qquad (11)$$

$$\Delta t_{C-H} \leq \frac{\Delta r^2}{4 + \frac{8\gamma*4}{\Delta r^2}} \qquad (12)$$

The above numerical stability criteria assume A-C and C-H equations with constant parameters. Varying parameters in the free energy density function (such as *T* and *h* in Eqn. 4) may influence the precise form of the stability condition of the explicit FD solver for solving C-H equation.

**3.2 Implicit FD solver**

This subsection presents the proposed implicit FD solver for A-C and C-H equations with varying parameters.

The unconditionally stable scheme must satisfy the following requirements: gradient stability (i.e., non-increasing total free energy over time) and von Neumann stability (*i.e.*, the errors made at one time step of the calculation will not cause the errors to be magnified in subsequent computations). Gradient stability is a subset of von Neumann stability for both A-C and C-H equations [40]. Thus, von Neumann stability is necessarily satisfied if gradient stability is satisfied. Accordingly, we only consider the gradient stability criterion when deriving the unconditionally stable implicit solver. We then verify unconditional numerical stability in selected numerical experiments.

The gradient stability condition has been well studied and summarized as Eyre's theorem [34], by following which we derive the unconditionally stable time integration method with varying $T$ and $h$, as detailed in the appendix. Assuming (without loss of generality) $M = 1$, we get

$$\phi_{t+\Delta t} + 2\Delta t(1-\xi)T\phi_{t+\Delta t} - \Delta t\gamma\nabla^2\phi_{t+\Delta t} = \phi_t - 2\xi T\Delta t\gamma\phi_t - 4\Delta t\phi_t^3 - h\Delta t \tag{13}$$

for the A-C equation, and

$$\phi_{t+\Delta t} - 2\Delta t\nabla^2\big((1-\xi)T\phi_{t+\Delta t}\big) + \Delta t\gamma\nabla^4\phi_{t+\Delta t} = \phi_t + \Delta t\nabla^2(2\xi T\phi_t + 4\phi_t^3 + h) \tag{14}$$

for the C-H equation. As derived in the appendix, gradient stability by Eyre's theorem requires

$$\begin{cases} \xi \geq \xi_{cri} = \frac{T(x,t)-12|\phi|_{max}^2}{2T(x,t)} & \text{if } T(x,t) < 0 \\ \xi \leq \xi_{cri} = -\frac{6|\phi|_{max}^2}{T(x,t)} & \text{if } T(x,t) > 0 \end{cases} \tag{15}$$

In other words, $h$ does not influence the stability of the numerical solver. Since the parameter $T$ is spatially and temporarily dependent, $\xi_{cri}$ will be as well.

*Direct solution to the implicit A-C scheme.*

The implicit scheme for the A-C equation in equation (13) can be rewritten in matrix format

$$(1+\boldsymbol{\epsilon})\odot\boldsymbol{\phi}_{t+\Delta t} + \eta\nabla^2\boldsymbol{\phi}_{t+\Delta t} = \boldsymbol{b} \tag{16}$$

with

$$\boldsymbol{\epsilon} = 2(1-\boldsymbol{\xi})\odot\boldsymbol{T}\Delta t, \quad \eta = -\Delta t\gamma, \quad \boldsymbol{b} = \boldsymbol{\phi}_t - 2\Delta t\boldsymbol{\xi}\odot\boldsymbol{T}\odot\boldsymbol{\phi}_t - 4\Delta t\boldsymbol{\phi}_t^3 - \boldsymbol{h}\Delta t \tag{17}$$

where $\odot$ denotes the Hadamard product or element-wise product of two matrices.

Representing the 2D order parameter field in the whole domain by a matrix $\boldsymbol{\Phi}$ with a dimension of $n \times m$, then, the Laplacian of $\boldsymbol{\Phi}$ in eq. (16) can be approximated with the finite-difference scheme as

$$\nabla^2\boldsymbol{\Phi} = \boldsymbol{A}_n\boldsymbol{\Phi} + \boldsymbol{\Phi}\boldsymbol{A}_m \tag{18}$$

where $\mathbf{A}_n$ and $\mathbf{A}_m$ are the matrices defined in Eqn. (7) or in Eqns. (8)-(10), accounting for BCs. $\mathbf{A}_n$ and $\mathbf{A}_m$ have the dimension of $n \times n$ and $m \times m$, respectively.

Since $\mathbf{A}_n$ and $\mathbf{A}_m$ are both symmetric matrices (i.e., $A_{ij} = A_{ji}$), we can express $\mathbf{A} = \mathbf{Q}\mathbf{D}\mathbf{Q}^T$, where $\mathbf{D}$ is a diagonal matrix whose elements along the diagonal are the eigenvalues of $\mathbf{A}$, $\mathbf{Q}$ is a matrix whose columns are the eigenvectors of $\mathbf{A}$ in accordance with the eigenvalues in $\mathbf{D}$, and $\mathbf{Q}^T$ is the transpose of $\mathbf{Q}$. Thus,

$$\nabla^2 \mathbf{\Phi} = \mathbf{Q}_n \mathbf{D}_n \mathbf{Q}_n^T \mathbf{\Phi} + \mathbf{\Phi} \mathbf{Q}_m \mathbf{D}_m \mathbf{Q}_m^T \tag{19}$$

Substituting into equation (16) yields

$$(1 + \boldsymbol{\epsilon}) \odot \mathbf{\Phi} + \eta(\mathbf{Q}_n \mathbf{D}_n \mathbf{Q}_n^T \mathbf{\Phi} + \mathbf{\Phi} \mathbf{Q}_m \mathbf{D}_m \mathbf{Q}_m^T) = \boldsymbol{b} \tag{20}$$

Left-multiplying both sides by $\mathbf{Q}_n^T$ and right-multiplying by $\mathbf{Q}_m$, we get

$$(1 + \boldsymbol{\epsilon}) \odot \mathbf{Q}_n^T \mathbf{\Phi} \mathbf{Q}_m + \eta(\mathbf{D}_n \mathbf{Q}_n^T \mathbf{\Phi} \mathbf{Q}_m + \mathbf{Q}_n^T \mathbf{\Phi} \mathbf{Q}_m \mathbf{D}_m) = \mathbf{Q}_n^T \boldsymbol{b} \mathbf{Q}_m \tag{21}$$

which may be rewritten as

$$(1 + \boldsymbol{\epsilon}) \odot \mathbf{Y} + \eta(\mathbf{D}_n \mathbf{Y} + \mathbf{Y} \mathbf{D}_m) = \mathbf{B} \tag{22}$$

with

$$\mathbf{Y} = \mathbf{Q}_n^T \mathbf{\Phi} \mathbf{Q}_m, \qquad \mathbf{B} = \mathbf{Q}_n^T \boldsymbol{b} \mathbf{Q}_m \tag{23}$$

Since $\mathbf{D}_n$ and $\mathbf{D}_m$ are both diagonal matrices, each element of $\mathbf{D}_n \mathbf{Y} + \mathbf{Y} \mathbf{D}_m$ in the $i$th row and $j$th column can be calculated as $\left(d_{n_i} + d_{m_j}\right) Y_{ij}$, where $d_{n_i}$ is the $i$th diagonal element of $\mathbf{D}_n$ and $d_{m_j}$ is the $j$th diagonal element of $\mathbf{D}_m$. Consequently, the left hand side of equation (22) can be written as $\mathbf{\Omega} \odot \mathbf{Y}$ with the elements of $\mathbf{\Omega}$ being the coefficients before each element $Y_{ij}$ of $\mathbf{Y}$. The elements of $\mathbf{\Omega}$ have the form

$$\Omega_{ij} = 1 + \epsilon_{ij} + \eta_{ij}(d_{n_i} + d_{m_j}) \tag{24}$$

With $\mathbf{\Omega}$ and $\mathbf{B}$ known, $\mathbf{Y}$ can be computed as

$$\mathbf{Y} = \mathbf{B}./\mathbf{\Omega} \tag{25}$$

where ./ represents element-wise division, i.e., $Y_{ij} = \frac{B_{ij}}{\Omega_{ij}}$. Once $\mathbf{Y}$ is obtained, the order parameter field $\mathbf{\Phi}$ may be calculated as

$$\boldsymbol{\Phi} = \mathbf{Q}_n \mathbf{Y} \boldsymbol{Q}_m^T \tag{26}$$

*Direct solution to the implicit C-H scheme.*

Similarly, to solve the implicit C-H scheme, we write

$$\boldsymbol{\epsilon} \odot \boldsymbol{\phi}_{t+\Delta t} + \boldsymbol{\eta} \odot \nabla^2 \boldsymbol{\phi}_{t+\Delta t} + \varepsilon \nabla^4 \boldsymbol{\phi}_{t+\Delta t} = \boldsymbol{b} \tag{27}$$

with

$$\boldsymbol{\epsilon} = 1 - 2\Delta t \nabla^2 [(1-\boldsymbol{\xi}) \odot \boldsymbol{T}], \quad \boldsymbol{\eta} = -2\Delta t (1-\boldsymbol{\xi}) \odot \boldsymbol{T}, \quad \varepsilon = \gamma \Delta t, \quad \boldsymbol{b} = \boldsymbol{\phi}_t + \nabla^2 [4\boldsymbol{\phi}_t^3 + 2\boldsymbol{T} \odot \boldsymbol{\xi} \odot \boldsymbol{\phi}_t + \boldsymbol{h}] \tag{28}$$

The Laplacian term in equation (27) can be approximated with the operator in (19). The fourth-order derivative may be written as $\nabla^4 \boldsymbol{\phi} = \nabla^2(\nabla^2 \boldsymbol{\phi})$. Applying the relationship in equation (19) two times gives

$$\nabla^4 \boldsymbol{\phi} = \mathbf{Q}_n \mathbf{D}_n^2 \mathbf{Q}_n^T \boldsymbol{\phi} + 2 \mathbf{Q}_n \mathbf{D}_n \mathbf{Q}_n^T \boldsymbol{\phi} \mathbf{Q}_m \mathbf{D}_m \mathbf{Q}_m^T + \boldsymbol{\phi} \mathbf{Q}_m \mathbf{D}_m^2 \mathbf{Q}_m^T \tag{29}$$

Inserting equations (19) and (29) into (27) and simplifying it gives

$$\boldsymbol{\epsilon} \odot \boldsymbol{\phi} + \boldsymbol{\eta} \odot (\mathbf{Q}_n \mathbf{D}_n \mathbf{Q}_n^T \boldsymbol{\phi} + \boldsymbol{\phi} \mathbf{Q}_m \mathbf{D}_m \mathbf{Q}_m^T) + \varepsilon (\mathbf{Q}_n \mathbf{D}_n^2 \mathbf{Q}_n^T \boldsymbol{\phi} + 2 \mathbf{Q}_n \mathbf{D}_n \mathbf{Q}_n^T \boldsymbol{\phi} \mathbf{Q}_m \mathbf{D}_m \mathbf{Q}_m^T + \boldsymbol{\phi} \mathbf{Q}_m \mathbf{D}_m^2 \mathbf{Q}_m^T) = \mathbf{b} \tag{30}$$

Left-multiplying $\mathbf{Q}_n^T$ and right-multiplying $\mathbf{Q}_m$ each term in equation (30), we get

$$\boldsymbol{\epsilon} \odot \mathbf{Y} + \boldsymbol{\eta} \odot (\mathbf{D}_n \mathbf{Y} + \mathbf{Y} \mathbf{D}_m) + \varepsilon (\mathbf{D}_n^2 \mathbf{Y} + 2 \mathbf{D}_n \mathbf{Y} \mathbf{D}_m + \mathbf{Y} \mathbf{D}_m^2) = \mathbf{B} \tag{31}$$

where $\mathbf{Y}$ and $\mathbf{B}$ have the same forms as equation (23).

Consequently, the left hand of eq. (31) can be written as $\boldsymbol{\Omega} \odot \mathbf{Y}$, and the elements of the coefficient matrix $\boldsymbol{\Omega}$ have the form

$$\Omega_{ij} = \epsilon_{ij} + \eta_{ij}(d_{n_i} + d_{m_j}) + \varepsilon(d_{n_i}^2 + 2 d_{n_i} d_{m_j} + d_{m_j}^2) \tag{32}$$

Then, the concentration matrix $\boldsymbol{\phi}$ can be computed in the same way, as shown in equations (25)-(26) for the implicit A-C solver.

Overall, the implicit C-H solver has the same workflow as the implicit A-C solver, but with different elements of matrices **b** and $\boldsymbol{\Omega}$ as shown in equations (28) and (32), respectively.

*Time stepping scheme*

There are many well-established adaptive time stepping schemes [41]–[46] for implicit phase-field solvers to pursue the best computational efficiency. Adaptive time-stepping is especially crucial for evolution problems varying characteristic times scales, as is expected in our application. However, as this topic is beyond the scope of the current work, we adopt a constant time step for all the following numerical tests.

## 4. Validation and numerical performance

We examine the accuracy and computational efficiency of the proposed implicit solver on examples with constant parameters $T$ and $h$. We then investigate its stability when applied to problems with temporally and spatially varying $T$ and $h$.

### 4.1 Implicit A-C solver

We benchmark the A-C solver on the evolution of an initially sharp straight interface, as shown in Figure 3a. The domain is a unit square discretized on a 500×500 uniform mesh. We use constant $T = -2, h = 0$ and set $M = 1, \gamma = 0.01$, $\Delta t_{explicit} = 5 \times 10^{-4}$, $\Delta t_{implicit} = 1 \times 10^{-2}$. The numerical steady state solutions from the explicit and proposed implicit FD solvers are compared to each other at three time points and to the analytical steady solution $\phi_{theory} = \tanh\left(\frac{x-0.5}{2\sqrt{\gamma}}\right)$ [47], as shown in Figure 3c. The proposed implicit A-C solver is accurate even using a 20 times larger $\Delta t$ than the explicit solver.

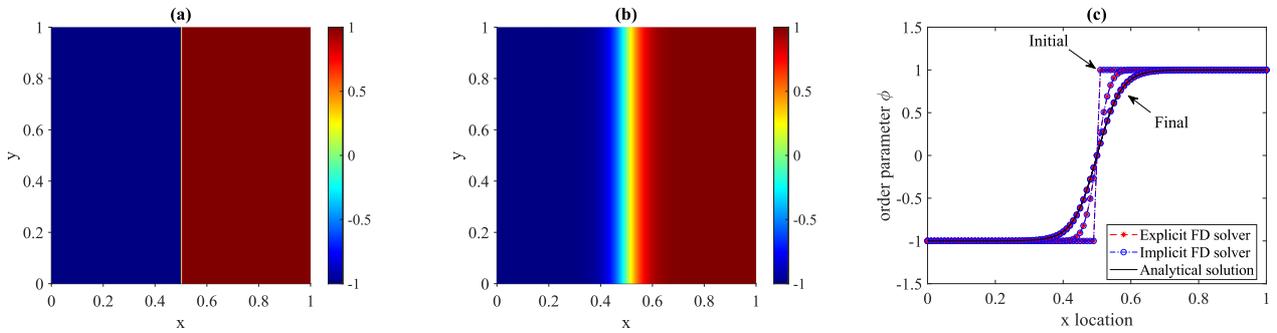

Figure 3 Benchmarking example of the smoothing of an initially sharp interface: (a) initial interface, (b) analytical solution for the final interface, (c) comparison of order parameter profiles along a line perpendicular to the interface obtained by the explicit solver, the proposed implicit solver, and analytical solution.

Next, we compare the numerical and gradient stability of the explicit solver and the proprosed implicit solver for different $\Delta t$. To investigate numerical stability, we introduced a random perturbation of magnitude $< 10^{-6}$ to every node of the initial state $\phi_0$ and examined how the $L^2$-norm of the resultant difference between the numerical solution with perturbations and a solution without random perturbations changes over time. The $L^2$-norm $\|\mathbf{\phi}_{M\times N} - \mathbf{\phi}'_{M\times N}\|_2$ is defined as $\sqrt{\frac{\sum_{i=1,M\ j=1,N}(\phi_{i,j}-\phi'_{i,j})^2}{M\times N}}$, where $\mathbf{\phi}$ and $\mathbf{\phi}'$ denote the numerical solutions of order parameter wihtout and with perturbation, respectively. If the difference norm decays or does not grow over time, the numerical solver is deemed stable numerically (and otherwise unstable). Gradient stability is deemed to be met if the total free energy, as defined in Eqn. (3), decays over time monotonously.

The stabiliy results of the explicit solver and proposed implicit solver are plotted in Figure 4 and Figure 5, respectively. According to the von-Neumann stability criterion in Eqn. (11), the explicit solver requires $\Delta t \leq 1.0 \times 10^{-3}$, which is evidenced by Figure 4. Whereas, as shown in Figure 5, the implicit solver is still stable both numerically and physically even when $\Delta t = 1000$. This outcome is consistent with the unconditional stability condition of the implicit solver, namely that the $\Delta t$ value of the implicit solver will no longer be limited by stability, but only by accuracy.

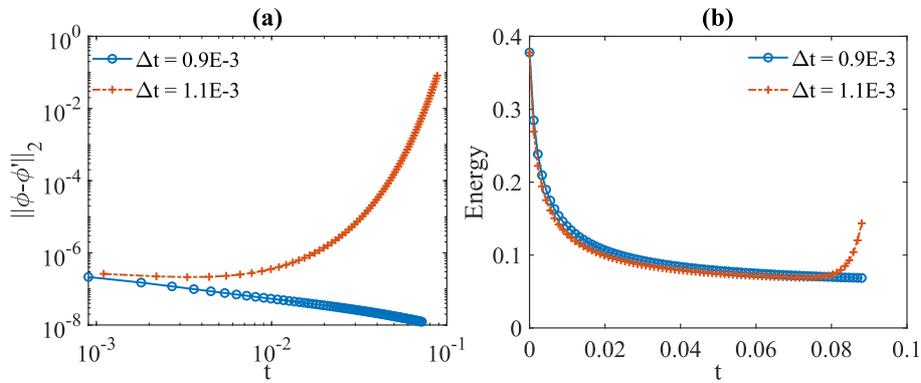

Figure 4 Stability of explicit A-C solver. (a) numerical stability (b) energy gradient stability.

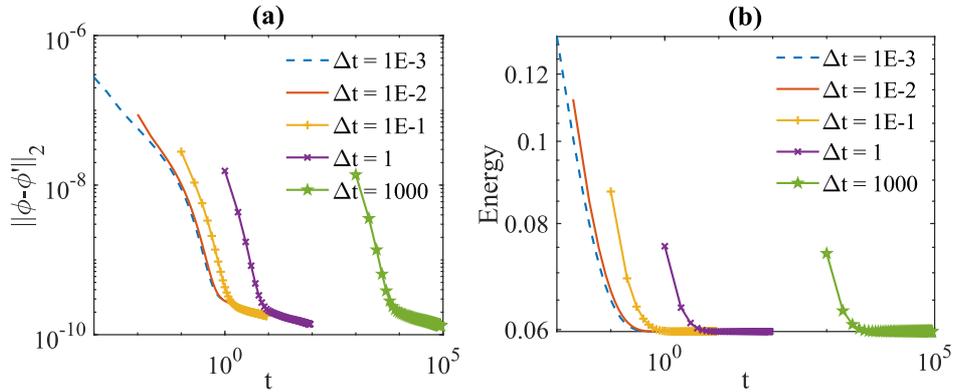

Figure 5 Stability of implicit A-C solver. (a) numerical stability, (b) energy gradient stability.

We then investigated the performance of the proposed implicit solver in computational efficiency relative to the explicit solver, as shown in Figure 6. As expected, the computation complexity of the explicit solver presented in section 3.1 is proportional to the number of degrees of freedom, DoF = $N \times N$, where $N$ is the number of grid points in $x$ and $y$ directions. Considering the necessary stability criterion in equation (11) applied to the explicit solver, $\frac{1}{\Delta t} \propto \frac{1}{\Delta r^2} \propto N^2 \propto DoF$, the explicit solver, therefore, owns a scalability of O(DoF²) overall, as evidenced by Figure 6. While for the implicit solver, its computation is dominated by multiplications of matrices with a dimension of $N \times N$, hence, the implicit solver has a computation complexty of $O(DoF^{1.5})$. Since $\Delta t$ in the implicit solver is constrainted by the accuracy requirement, which is proportional to $\Delta x$ or $DoF^{0.5}$, we conclude that the implicit solver owns the same $O(DoF^2)$ scalability as the explicit solver, as evidenced by Figure 6. Eventually, thanks to significantly larger allowable $\Delta t$, the implicit solver is ~7 times faster than the explcit solver for numerical models with large DoF.

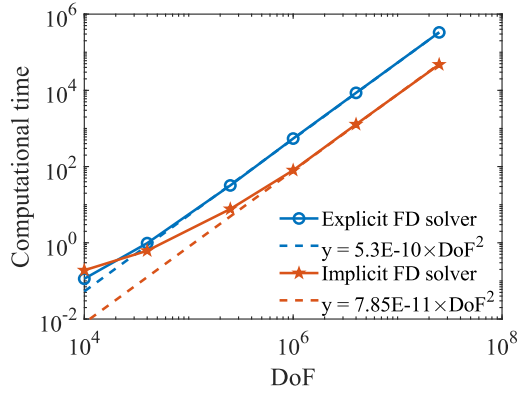

Figure 6 Computational efficiency comparison of explicit solver and implicit solver

## 4.2 Implicit C-H solver

This subsection validates the accuracy and stability condition of the proposed implicit C-H solver with two benchmarking examples, and then evaluates its computational efficiency performance relative to the explicit FD solver and the FEM solver integrated into the MOOSE software [48] developed by Idaho National Laboratory. These two examples different from the one for A-C solver's validation are used, because they are more suitable to show the conservation feature of $\phi$ in C-H equation.

The first example solves for the evolution of a square-shaped inclusion of one phase embedded in a matrix of another phase. Based on capillarity arguments, the initially square shape is expected to change gradually into a circle of identical area. The domain is discretized as $100 \times 100$. We used constant $T = -2, h = 0$ and set $M = 1, \gamma = 4 \times 10^{-4}$, $\Delta t_{implicit} = 1 \times 10^{-4}$. The numerical results in Figure 7 show that the implicit solver predicts an accurate interface profile while rigidly holding the conservation of concentration during evolution as needed.

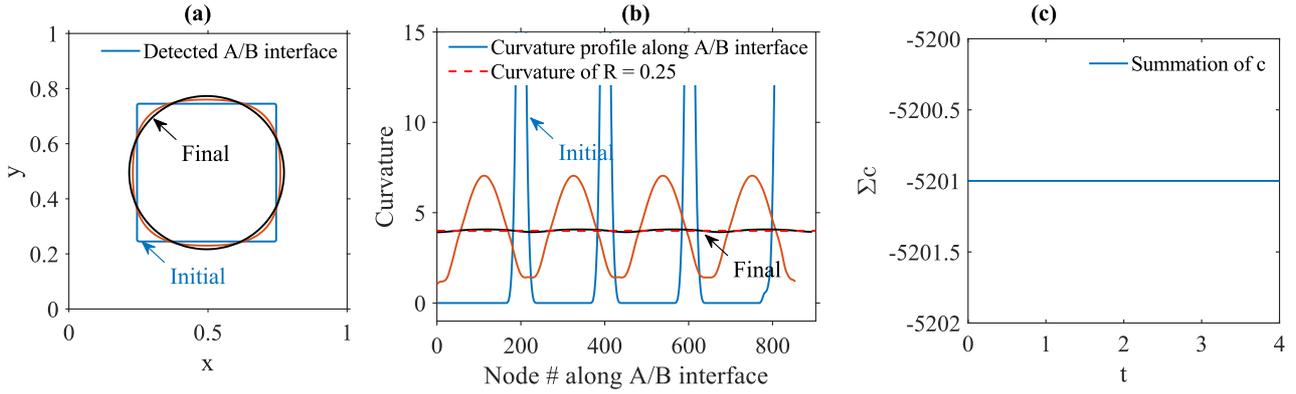

Figure 7 Evolution of a squared phase by implicit C-H solver. (a) interface profiles, (b) interface curvature [8], (c) time history of concentration integral in whole domain.

The stability testing results of the explicit solver and implicit solver are shown in Figure 8 and Figure 9, respectively. From Figure 8, we observe that the explicit solver is unstable when $\Delta t > 8 \times 10^{-7}$, which agrees with the criterion in equation **Error! Reference source not found.**. The implicit solver is stable unconditionally.

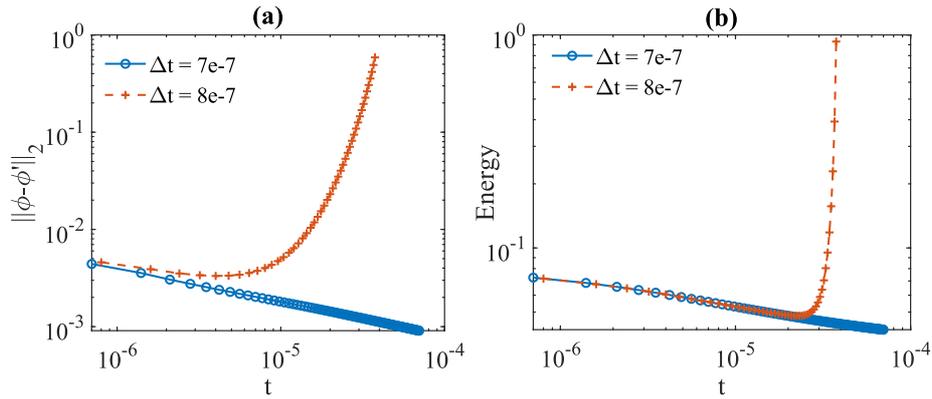

Figure 8 Stability of explicit C-H solver. (a) numerical stability, (b) energy stability.

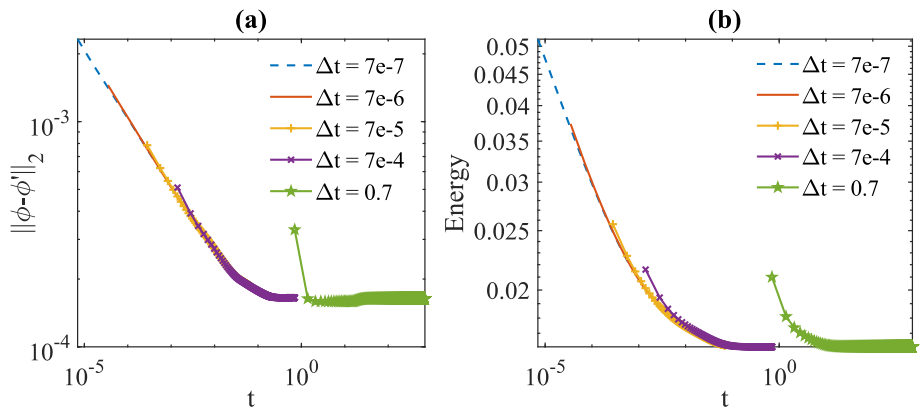

Figure 9 Stability of implicit C-H solver. (a) numerical stability, (b) energy stability.

For computation complexity, the implicit solver is O(DoF$^{1.5}$) complexity vs O(DoF) of explicit FD solver, as described in subsection 4.1. In time integration, the explicit one is limited by the numerical stability criterion which requires $1/\Delta t \propto 1/\Delta r^4 \propto DoF^2$. While the implicit solver is unconditional stable, and $\Delta t$ is only limited by the accuracy requirement which demands $1/\Delta t \propto 1/\Delta r \propto DoF^{0.5}$. This yields an advantageous scaling of O(DoF$^2$) for the implicit solver, compared to O(DoF$^3$) for the explicit solver. This conclusion agrees with the numerical testing results presented in Figure 10.

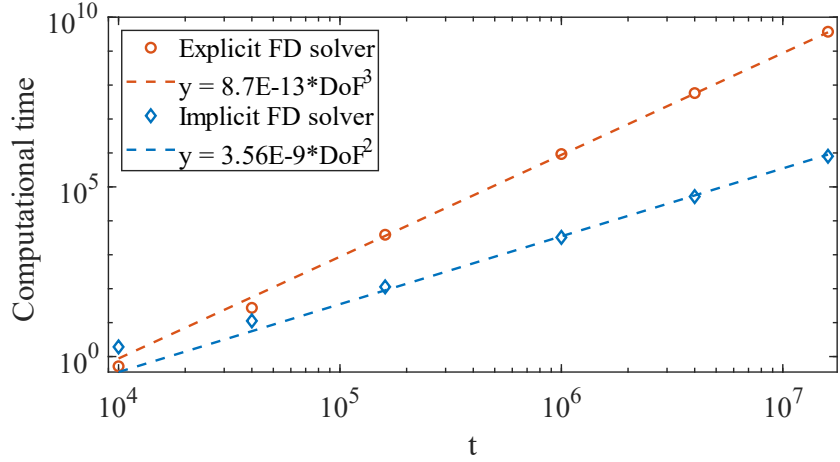

Figure 10 Computational efficiency comparison of explicit and implicit C-H solvers.

Finally, we compare the effectiveness and computational efficiency of the proposed implicit C-H solver to the well-established FEM-based platform, MOOSE [48]. *** Figure 11 shows that the solution obtained from the proposed implicit solver matches very well with that of the FEM solution by MOOSE. Under the same configuration of fixed $\Delta t$ and uniform grid, the implicit solver achieves a 5 times faster computation, compared to MOOSE. More detail about the geometry and MOOSE solver can be found from [8].

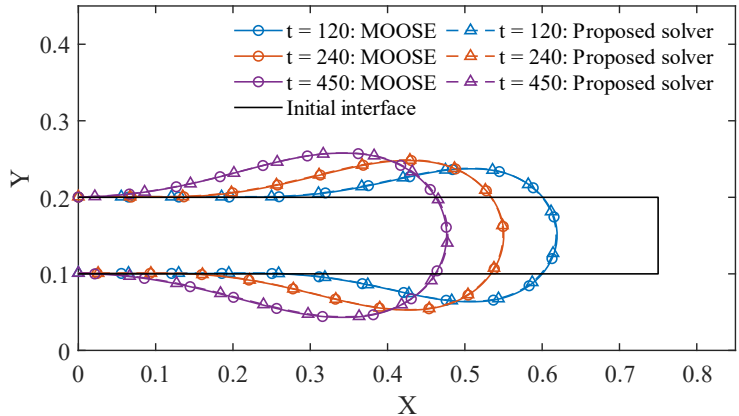

Figure 11 Results comparison of layer retraction by MOOSE and the proposed implicit solver. In simulation, a $400 \times 800$ grid is adopted with $M = 1, \gamma = 4E - 8, T = -2, h = 0, \Delta t = 1 \times 10^{-2}$.

**4.3 Examples with spatially and temporally varying T and h**

In this section, the stability condition of the proposed implicit solver is investigated when applied to the A-C and C-H equations with varying $T(x,t)$ and $h(x,t)$. The numerical model has a $50 \times 50$ grid with $M = 1, \gamma = 0.01$.

The first example uses the A-C equation to addresses the banded distribution presented in Figure 12a. Here, $T$ and $h$ have different values in the blue and red regions. However, within a given color (blue or red), $T$ is identical and $h$ as well. Both $T$ and $h$ change over time. The assigned values of $T$ and $h$ at different times and locations are listed in Table 1. Under the given configuration, the critical time step above which the explicit solver is unstable is $\Delta t_{cri} = 0.01$, based on the criterion in Eqn. (11). This is verified by Figure 13 showing that the explicit solver is stable numerically (in Figure 13e) and physically (Figure 13f) when $\Delta t < \Delta t_{cri}$, but unstable when $\Delta t \geq \Delta t_{cri}$. In contrast, the implicit one is still stable both numerically and physically as $\Delta t > \Delta t_{cri}$ without apparent accuracy loss.

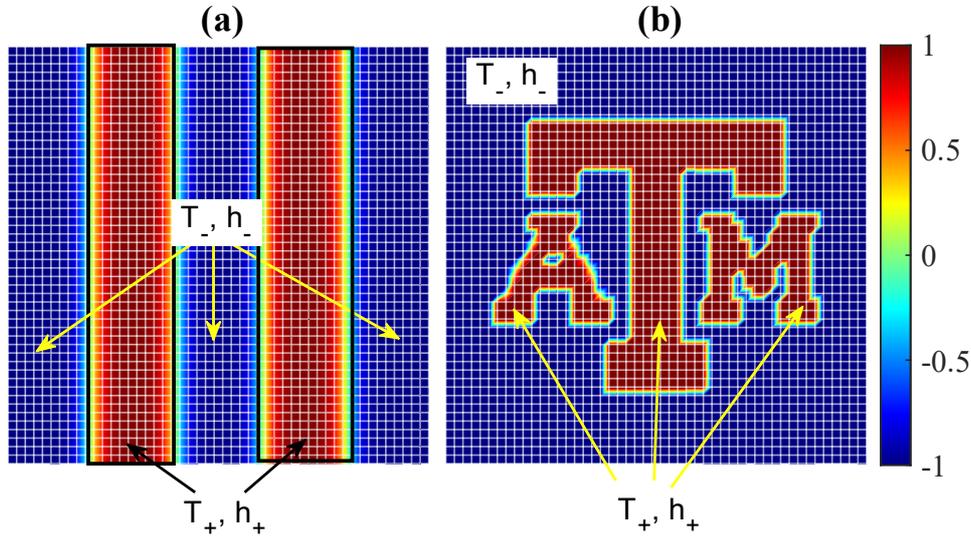

Figure 12 Configuration of $T(x,t)$ and $h(x,t)$ in the A-C example (a) and C-H example (b). $(T_-, h_-)$ is employed in the blue region, $(T_+, h_+)$ is employed in the red region. $T$ and $h$ change over time as listed in Table 1 and Table 2.

Table 1. List of $T$ and $h$ values within four time periods in the A-C example. $(T_-, h_-)$ is employed in the blue region, $(T_+, h_+)$ is employed in the red region.

| Time periods | $T_+$ | $T_-$ | $h_+$ | $h_-$ |
|---|---|---|---|---|
| $t \in [0.00, 0.05]$ | 0.49 | -4.13 | -11.76 | 3.24 |
| $t \in (0.05, 0.10]$ | 4.18 | 0.81 | -5.15 | 2.41 |
| $t \in (0.10, 0.15]$ | 1.91 | -2.51 | -9.85 | 4.88 |
| $t \in (0.15, 0.20]$ | 2.75 | 2.23 | -15.23 | 8.05 |

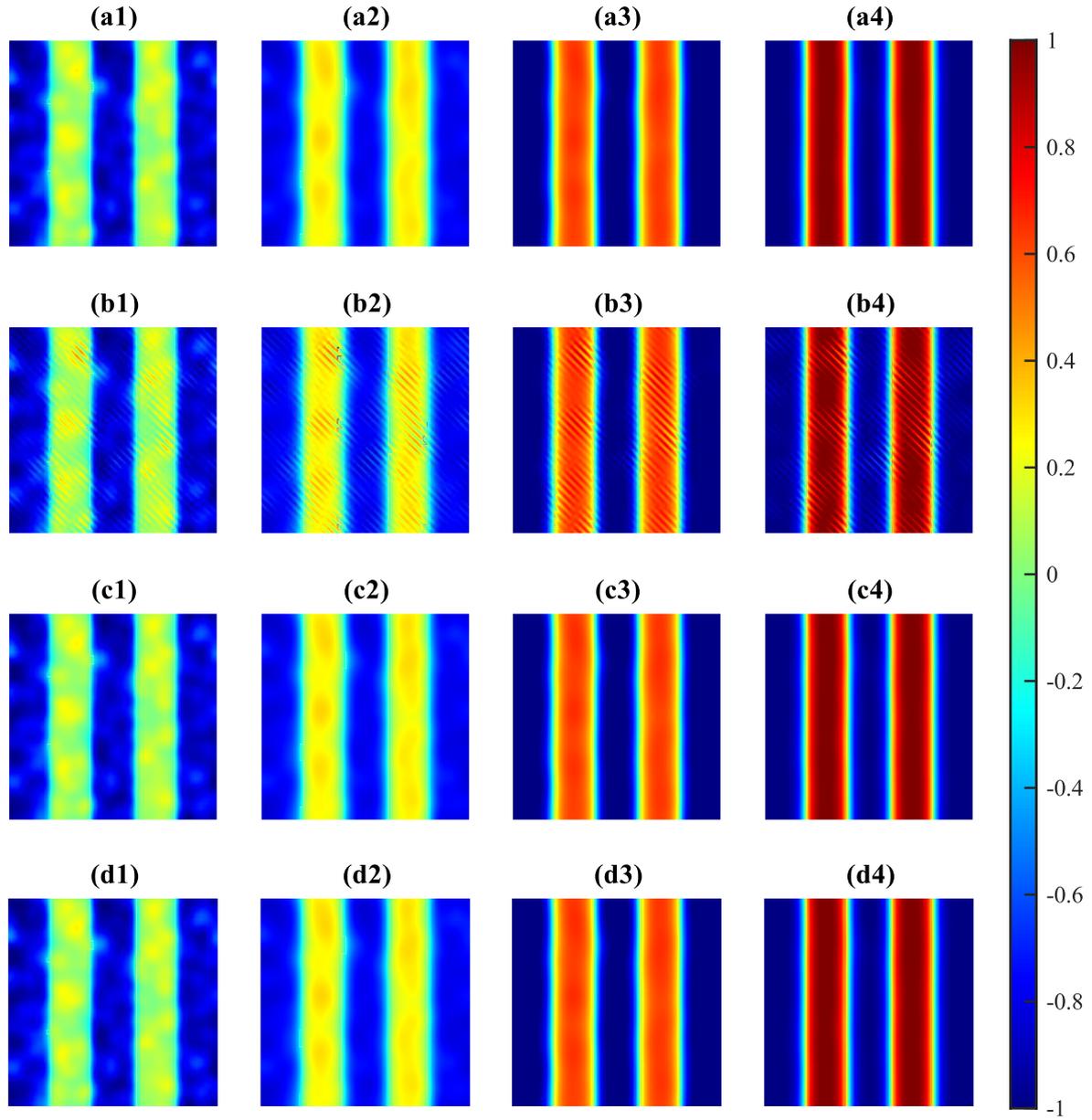

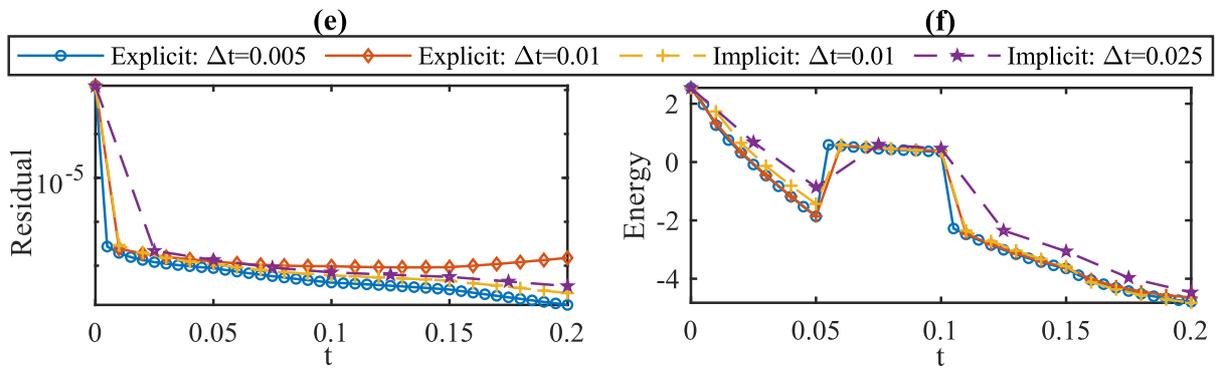

Figure 13 Results of A-C simulations by the explicit solver with (a) $\Delta t = 0.5\Delta t_{cri}$, (b) $\Delta t = \Delta t_{cri}$, and by the implicit solver with (c) $\Delta t = \Delta t_{cri}$, (d) $\Delta t = 2.5\Delta t_{cri}$. (e) plots the change of residual caused by initial perturbation over t while (f) plots the change of total energy over t.

The stability of implicit C-H solver is evaluated with an example where the temporally varying $T$ and $h$ have a more complex distribution, as presented in Figure 12b and listed in Table 2. The numerical results presented in Figure 14 show that the explicit solver is unstable when $\Delta t = 4 \times 10^{-7}$ under the given varying $T$ and $h$. By contrast, the implicit solver is stable even when the time step is as large as $1 \times 10^{-4}$ while maintaining order parameter conservation.

Table 2. Table 3 List of $T$ and $h$ values within four time periods in the A-C example. $(T_-, h_-)$ is employed in the blue region, $(T_+, h_+)$ is employed in the red region.

| Time periods | $T_+$ | $T_-$ | $h_+$ | $h_-$ |
| --- | --- | --- | --- | --- |
| $t \in [0.00,0.05]$ | 9.44 | 4.72 | -22.84 | 13.41 |
| $t \in (0.05,0.10]$ | -7.52 | 2.51 | -19.08 | -1.03 |
| $t \in (0.10,0.15]$ | 9.73 | 12.65 | -23.46 | 29.30 |
| $t \in (0.15,0.20]$ | 25.31 | 25.37 | -54.69 | 54.57 |

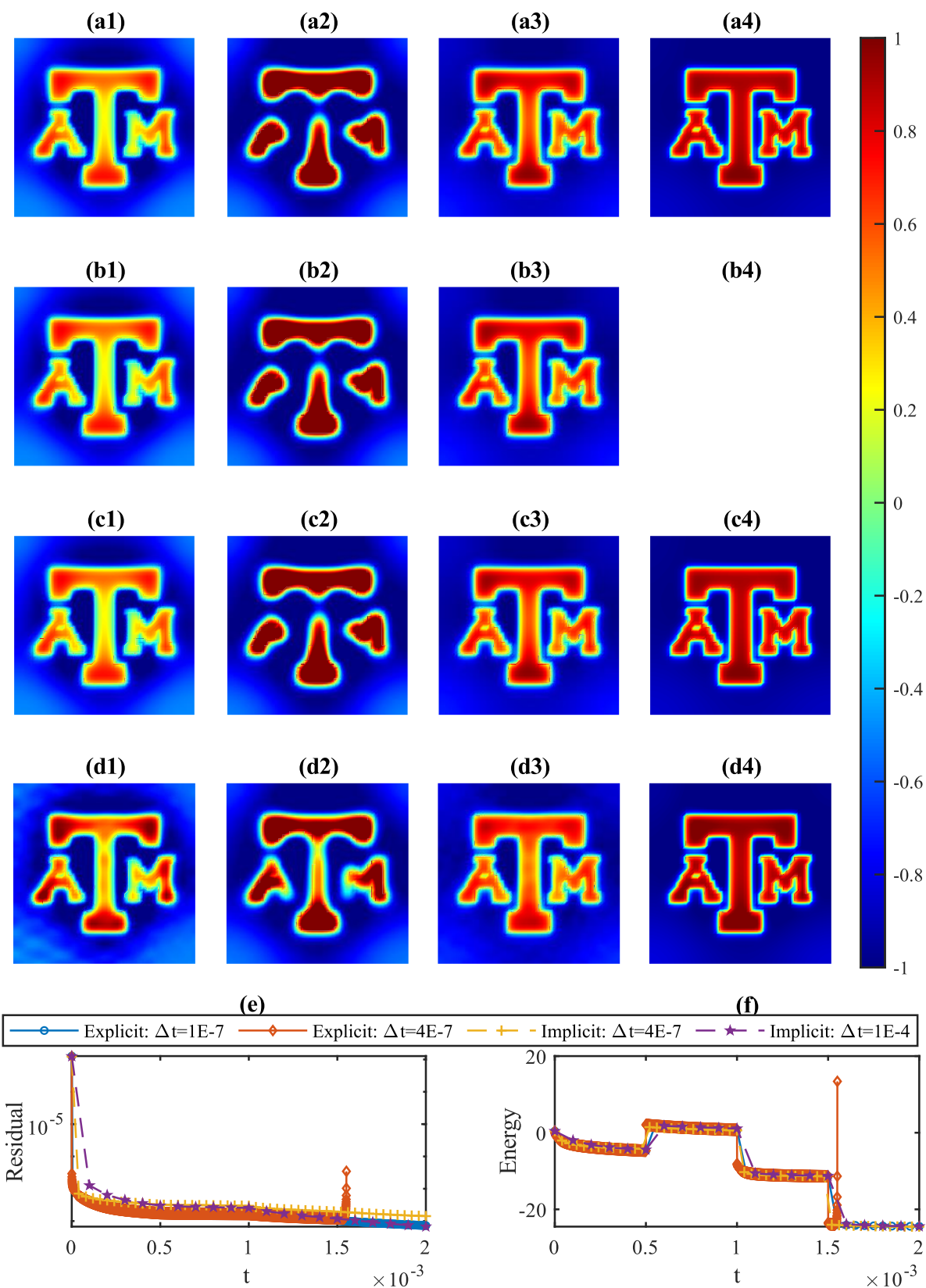

Figure 14 Results of C-H simulations by the explicit solver with (a) $\Delta t = 1 \times 10^{-7}$, (b) $\Delta t = 4 \times 10^{-7}$, and by the implicit solver with (c) $\Delta t = 4 \times 10^{-7}$, (d) $\Delta t = 1 \times 10^{-4}$. (e) plots the change of residual over t while (f) plots the change of total energy over t.

Additional numerical experiments (not detailed here) demonstrate that the implicit solvers are still stable even at $\Delta t = 100{,}000 \Delta t_{cri}$, though may not be accurate.

## 5. Discussion and Conclusions

We have proposed implicit A-C and C-H solvers that are accurate, unconditionally stable, computationally efficient, and applicable to problems with spatially and temporally varying parameters. We have demonstrated these properties using both theoretical analysis and numerical experiments.

Compared to explicit FD solvers, the proposed implicit solver is unconditionally stable, which allows it to use longer time steps. Its efficiency advantage is especially significant for the fourth-order Cahn-Hilliard models with large numbers of DoFs. Indeed, the implicit C-H solver scales as O(DoF$^2$) while the explicit solver scales as O(DoF$^3$). Moreover, some implicit solvers use iterative numerical approaches to solve the FD equations, so their computational efficiency depends on the convergence rate of the employed iterative method and memory storage requirements. By contrast, the implicit solver proposed here obtains a direct solution, further improving its efficiency and predictability. This is why the proposed implicit solver outperforms the iterative FEM-based MOOSE framework in computation efficiency.

Compare to the FEM-based MOOSE, the FD-based implicit solver is much easier in boundary treatments. This is because boundary conditions are always expressed as derivatives of independent variables, while FD solvers solve the differential form of governing equations directly. So, in the boundary treatment aspect, our implicit FD solver has an advantage over FEM, which solves the integral form of governing equations and therefore has challenges in boundary treatment.

One drawback of our implicit solvers is that they are presently only applicable to phase-field problems with constant mobility. For the problems with order parameter-independent mobility, an implicit solver may also be developed following the deduction in appendix. However, a direct solution of the corresponding discrete scheme might not be achievable. Furthermore, direct solutions for the proposed implicit solver are only achievable for problems with a periodic BC, Neumann BC, or PBC/NBC mixed boundary meaning some directions are periodical while some other directions are NBC. This is because the calculation step in equation (19) requires a "definite positive and symmetric" coefficient matrix **A** defined in equations (8)-(10). Nevertheless, considering that PBCs are the most widely used boundary conditions in phase-field modeling, the direct solver is sufficient for many applications.

Finally, all the results presented here used uniform meshes. In practice, many solvers (e.g., in MOOSE) achieve improved performance by using non-uniform, adaptive meshes. These techniques are of great benefit in phase field models, which demand high resolution at interfaces and not necessarily elsewhere. To adapt our implicit solver to non-uniform meshes, we may turn to the Gradient Smoothing Method [49], though ease of implementation of our FD solvers may suffer, as a result.

# Appendix

This section derives the unconditionally gradient stable scheme presented in (13)-(15).

The free-energy functional is defined as

$$F \equiv \int dV [\tfrac{1}{2}\gamma(\nabla\boldsymbol{\phi})^2 + (\boldsymbol{\phi}^4 + \boldsymbol{T}\boldsymbol{\phi}^2 + \boldsymbol{h}\boldsymbol{\phi})] \tag{33}$$

Eyre's theorem [34] says that an algorithm will be *unconditionally gradient stable*, for both A-C and C-H equations, if the free energy functional $F$ can be split appropriately into the so-called 'contractive' and 'expansive' parts, $F = F^C + F^E$, and the contractive parts are treated implicitly and the expansive parts are treated explicitly. The necessary condition on the splitting is identical for A-C and C-H equations and may be stated by introducing the Hessian matrices,

$$M_{ij} = \frac{\partial^2 F}{\partial\phi_i \partial\phi_j}, \quad M_{ij}^E = \frac{\partial^2 F^E}{\partial\phi_i \partial\phi_j}, \quad M_{ij} = \frac{\partial^2 F^C}{\partial\phi_i \partial\phi_j} \tag{34}$$

where $i, j$ denote lattice sites.

According to the Eyre's theorem [34], the unconditional gradient stability requires that i) all eigenvalues of $\mathbf{M}^E$ must be nonpositive, and all eigenvalues of $\mathbf{M}^C$ must be non-negative, ii) the largest eigenvalue, $\lambda_{max}^E$, of $\mathbf{M}^E$ must be not greater than a half of the smallest eigenvalue, $\lambda_{min}$, of $\mathbf{M} = \mathbf{M}^E + \mathbf{M}^C$, i.e.,

$$\lambda_{max}^E \leq \tfrac{1}{2}\lambda_{min} \tag{35}$$

To identify the proper splitting, we break the free energy $F$ into four parts

$$F^{(1)} = \sum_i T\phi_i^2, \quad F^{(2)} = \sum_i \tfrac{1}{2}\gamma(\nabla\phi_i)^2, \quad F^{(3)} = \sum_i \phi_i^4, \quad F^{(4)} = \sum_i h\phi_i \tag{36}$$

First, $M_{ij}^{(1)} = 2T\delta_{ij}$ and $M_{ij}^{(3)} = 12\phi_i^2 \delta_{ij}$, where $\delta_{ij}$ is the Kronecker function, have eigenvalues equal to the diagonal elements. Next, $M_{ij}^{(2)} = (-\gamma\nabla^2)_{ij}$ can always be diagonalized by going to Fourier space, which yields that the eigenvalues of $\mathbf{M}^{(2)}$ are strictly non-negative when $\gamma$ is positive. Finally, $M_{ij}^{(4)} = 0$. We parametrize the splitting via

$$F^E = \sum_{k=1}^{4} a_k F^{(k)}, \quad F^C = \sum_{k=1}^{4} (1 - a_k) F^{(k)} \tag{37}$$

Now, the values of $a_i$ in above equation must be chosen to satisfy the two stability requirements (i-ii). If we take $a_2 = 0$, $a_3 = 1$, and $a_4 = 1$ for simplicity, yielding $F^E = a_1 F^{(1)} + F^{(3)} + F^{(4)}$, then the eigenvalues of $\mathbf{M}^E$ are of the form $2a_1 T + 12\phi_i^2$ which should be non-positive (by gradient stability requirement (i)) and its greatest eigenvalue, i.e., $\lambda_{max}^E = 2a_1 T + 12|\phi|_{max}^2$, cannot be greater than $1/2\,\lambda_{min}$ (by requirement (ii)).

Considering the eigenvalues of the four parts of free energy term, $\lambda_{min}$ depends on the value of $T$. If $T < 0$, then $\lambda_{min} = 2T < 0$, and the stability requirements are met by

$$a_1 \geq \frac{T - 12|\phi|_{max}^2}{2T} \tag{38}$$

The above condition satisfies requirement (i).

Otherwise, if $T > 0$, then $\lambda_{min} = 2T > 0$. In this case, $\lambda_{max}^E$ should be not greater than 0 because of requirement (i), which gives

$$a_1 \leq -\frac{6|\phi|_{max}^2}{T} \tag{39}$$

Consequently, Eyre's theorem guarantees that the numerical scheme will be unconditionally gradient stable by treating $F^E$ explicitly and $F^C$ implicitly with the conditions in (38) and (39).

# Acknowledgements

This work was supported by the NSF CDS&E program under grant number 1802867. Computational resources were provided by the High-Performance Research Computing (HPRC) center at Texas A&M University.